\documentclass[11pt]{article}   
 \pdfoutput=1 
 \usepackage{graphicx}
\usepackage{amsthm,amsmath,amssymb}   
\usepackage{color}
\usepackage{url}
\usepackage{enumerate,enumitem}
\usepackage{epsfig,subfigure}
\newtheorem{theorem}{Theorem}[section]

\newtheorem{proposition}[theorem]{Proposition}

\theoremstyle{remark}

\title{\LARGE A half-automated study \\ of a 2-parameter family of integrals}

\date{\today}

\begin{document}

\maketitle
\begin{center}
{\bf David G. Zeitoun}\\[1mm]
Talpiot College of Education \\
Holon, Israel \\ 
{\tt ed.technologie@gmail.com}\\[1mm]
{\bf  Thierry Dana-Picard}\\[1mm]
Jerusalem College of Technology and Jerusalem Michlala College \\ 
Jerusalem, Israel \\
{\tt  ndp@jct.ac.il}
\end{center}

\vspace*{1mm}

\paragraph{Abstract:}
\small
{\it
The study of some parametric integrals is presented with a combined approach of analytical development, the usage of a Computed Algebra System (CAS) and of the Online Encyclopedia of Integer Sequences. 
The methodology for the solution includes a) an analytical investigation for the study of the parametric integral, b) computations with a CAS of the integral for specific values of the parameter, c) investigation of the connection between the integral and special functions or classical numbers, and d) derivation of a general algorithm for the complete computation of the parametric integral. The central example of the paper is the parametric integral 
\begin{equation*}
\label{the general integral}
I_{n}^{(p)}=\int_0^{\pi /4} x^p \tan^n x\; dx,
\end{equation*}
The work reveals a connection of  this parametric integral with Catalan numbers.}
\normalsize
\vspace*{1mm}

\section{Introduction.}
\label{intro}
The computation of parametric definite integrals is an interesting mathematical field (e.g. see \cite{moll1,moll2}). On the one hand, parametric integrals help students to grasp more abstract situations than in the cases without parameters. On the other hand, in domains such as Physics and Engineering, numerous phenomenons have a mathematical translation into such integrals; for example, see \cite{soil}.  Moreover, such computations may lead to formula for indefinite series, and also to derive combinatorial identities and integral representations for combinatorial objects; for example, see \cite{DP-Catalan,stirling,bernstein book,wallis-catalan,Qi et al 2017,Qi Elemente 2017,Yin and Qi} and the papers in reference there. In the specific case of Catalan numbers, different presentations and applications are given in \cite{koshy}.

Computing the integrals by hand may be unilluminating. The usage of technology has been analyzed and discussed for a long time, not only for the domain that we are interested in here.
Buchberger \cite{buchberger} says that a Computer Algebra System (CAS) should be used only where the user knows how to perform the computations by hand, not as a black box where the user has no knowledge on what are the processes at work inside. A good example of this approach is in the book \cite{adams and loustaunau} where the described algorithms  are first "run by hand", and then applied with their version implemented in software\footnote{The book is devoted to Gr\"obner bases and their applications. This is not the topic of our work here, but the didactic methodology is similar to what we claim}. With a black box usage, the CAS is a facilitator to obtain a result, not to develop more mathematical knowledge. For a short presentation of this \emph{White-Box / Black-Box Principle}, we refer to \cite{wbbb}. 

Nevertheless, but as already seen in \cite{bypass}, it happens that the CAS is used to bypass a lack of theoretical knowledge, then proceeding backwards the gap can be filled. In \cite{drijvers}, Drijvers  notes that generally a CAS works as a black box, it does not provide insight into the algorithms. Advanced versions of a CAS may provide some insight, either as a step-by-step command, as in the Derive software (in version 6.1, released a long time ago, which was the last release), or via tutorials such as Maple's tutorials. These tutorials work in an interactive way, providing hints for every single step, sometimes more than one. The user then chooses and checks the efficiency of the proposed step. A suitable usage of such a tutorial provides a good scaffolding to the user, which has to be gradually faded. Note that the same tutorials have a shortcut giving the final answer immediately. In this case, the CAS has been used as a black box. Anyway, this usage is highly personal , and each user will develop his/her own instrumental genesis \cite{artigue,trouche}. A good balance between paper-and-pencil and CAS assisted work is crucial.
 
In this paper, we study a 2-parameter family of definite integrals. For small values of the parameters, computations may be performed by hand, but this is time-consuming. CAS can be used after a few computations, without contradicting Buchberger's point of view. This may be a partial solution only, as we will see in the example developed in Section \ref{section induction}. In certain cases, the output of an automated computation shows an integral which cannot be computed analytically. From that point, either numerical methods should be used, or an induction formula has to be looked for, as in \cite{stirling,wallis-catalan} and other cases referenced there. 

The numerical results obtained with a CAS have to be exploited. The investigation made a strong usage of the Online Encyclopedia of Integer Sequences\footnote{\url{http://www.oeis.org}.}. This proved a really efficient approach, enabling to arrive to conjectures, which have then to be proven. 

We study here a generalization of a parametric integral described in \cite{wallis-catalan}, namely the 2-parameter integral defined by:
\begin{equation}
\label{the general integral}
I_{n}^{(p)}=\int_0^{\pi /4} x^p \tan^n x\; dx,
\end{equation}
where the parameters $p$ and $n$ are non negative integers.
First we derive by hand an induction formula, then we try to have more insight into the 2-parameter sequence of integrals. For this we use a Computer Algebra System, here the Maple package and its tutorials, and the Online Encyclopedia of Integer Sequences.

Automated methods have been analyzed and documented for questions in Geometry, Analytic Geometry and in Differential Geometry (for example, see the survey \cite{DP-automated}). The important novelty here is that, for parametric integrals no visual intuition, actually almost no intuition at all, can be afforded. We should mention the example of $\int x \ln x \; dx$ for which the choice of the functions for an integration by parts is counter-intuitive.  The present study requests more abstract thinking. After all, parametric integrals are already a more abstract object than "ordinary" definite integrals.As we will see, "the computer offers to provide scaffolding both to enhance mathematical reasoning and to restrain mathematical error" (J. Borwein in\cite{borwein}).

\section{An induction formula}
\label{section induction}

For integrals where the integrand contains a power of the tangent function, it is natural to look for an induction with an increment of 2 in the power.
We have:
\begin{align*}
I^{(p)}_n+I^{(p)}_{n-2} & = \int_0^{\pi /4} x^p \tan^n x \; dx+\int_0^{\pi /4} x^p \tan^{n-2} x \; dx\\
\quad &= \int_0^{\pi /4} x^p ( \tan^n x + \tan^{n-2}x )\; dx\\
\quad & = \int_0^{\pi /4} x^p ( \tan^2 x + 1 ) \; \tan^{n-2}x ; dx\\
\end{align*}
Integration by parts of the last expression, using the functions:
\begin{equation}
f(x) = x^p  \; \text{and} \; g(x)= \frac{\tan^{n-1}x}{n-1} 
\end{equation}
whence
\begin{equation}
f'(x) = p x^{p-1} \; \text{and} \; g'(x) = \tan^{n-2}x ( \tan^2 x + 1 )
\end{equation}
This leads to the following recurrence relation:
\begin{equation}
\label{fundamentalnp}
I^{(p)}_n+I^{(p)}_{n-2}= \frac{1}{n-1} \left( \frac{\pi}{4} \right)^p +\frac{p}{n-1}I^{(p-1)}_{n-1}, \; p \geq  1; \;  n\geq 1.
\end{equation}
We can also derive general  recurrence formulas for both the sum $\underset{n}{\sum} I^{(p)}_n$ for a given $p$, and for $I^{(p)}_n$.
\begin{proposition}
When $n$ is even, we denote $n=2l$,  and the following holds:
\begin{equation}
\label{induction 2l}
I^{(1)}_{2l} = \frac{\pi}{4} \underset{k=1}{\overset{2l-1}{\sum}} \frac {(-1)^k}{k}
+  \underset{k=1}{\overset{2l-1}{\sum}} \frac {(-1)^k I_{k}}{k}- \frac{1}{2} \ \left( \frac{\pi}{4} \right)^2
\end{equation}
\end{proposition}

\subsection{The case of even $n$}

\subsubsection{Analytic work}

For a given integer $p \geq \ 1$ and if $n$ is even:

\begin{align*}
I^{(p)}_n + 2 (I^{(p)}_{n-2}+I^{(p)}_{n-4}+ \dots + I^{(p)}_2 ) +I^{(p)}_0 & = \left( \frac{\pi}{4} \right)^p \underset{k=1}{\overset{n-1}{\sum}} \frac {1}{k}
+ p \underset{k=1}{\overset{n-1}{\sum}} \frac {I^{(p-1)}_{k}}{k}\\
\quad & = ( \Psi (n) + \gamma)\left(\frac{\pi}{4} \right)^p+p \underset{k=1}{\overset{n-1}{\sum}} \frac {I^{(p-1)}_{k}}{k}.
\end{align*}
By addition and substraction of the iterative formula in Equation (\ref{fundamentalnp}), we derive the following formula:
\begin{equation}
I^{(p)}_n + I^{(p)}_0 = \left(\frac{\pi}{4} \right)^p \underset{k=1}{\overset{n-1}{\sum}} \frac {(-1)^k}{k}
+ p \underset{k=1}{\overset{n-1}{\sum}} \frac {(-1)^k I^{(p-1)}_{k}}{k}.
\end{equation}
The first integral in the sequence is
\begin{equation}
 I^{(p)}_0 = \frac{1}{p+1} \ \left(\frac{\pi}{4} \right)^{(p+1)},
\end{equation}
therefore, the following holds:
\begin{equation}
\label{eq 1}
I^{(p)}_n = \left( \frac{\pi}{4} \right)^p \underset{k=1}{\overset{n-1}{\sum}} \frac {(-1)^k}{k}
+ p \underset{k=1}{\overset{n-1}{\sum}} \frac {(-1)^k I^{(p-1)}_{k}}{k}- \frac{1}{p+1} \ \left( \frac{\pi}{4} \right)^{(p+1)}
\end{equation}
Equation (\ref{eq 1}) means that $I_n^{(p)}$ is a linear function of the previous integrals $I_n^{(k)}$, for $k<p$.

If $n$ is even and  $p=0$, we obtain the following identity:
\begin{equation}
I^{(0)}_n + I^{(0)}_0 = \underset{k=1}{\overset{n-1}{\sum}} \frac {I_{k}}{k}.
\end{equation}
Using the previous derivation, we obtain $ I^{(0)}_0= \frac{\pi}{4}$, whence:
\begin{equation}
I^{(0)}_n  = \underset{k=1}{\overset{n-1}{\sum}} \frac {I_{k}}{k} - \frac{\pi}{4}.
\end{equation}

If $n$ is even and  $p=1$, we obtain:
\begin{equation}
I^{(1)}_n + 2 (I^{(1)}_{n-2}+I^{(1)}_{n-4}+ \dots + I^{(1)}_2 ) +I^{(1)}_0 = \frac{\pi}{4}( \Psi (n) + \gamma)
+ \underset{k=1}{\overset{n-1}{\sum}} \frac {I_{k}}{k},
\end{equation}
where  $\Psi (n)$ is the  \emph{Digamma function}\footnote{\url{https://www.maplesoft.com/support/help/Maple/view.aspx?path=Psi}} and $\gamma$ is the Gamma number.

Also a general expression of $I^{1}_n$ may be derived:
\begin{equation}
\label{induction n even}
I^{(1)}_n = \frac{\pi}{4} \underset{k=1}{\overset{n-1}{\sum}} \frac {(-1)^k}{k}
+  \underset{k=1}{\overset{n-1}{\sum}} \frac {(-1)^k I_{k}}{k}- \frac{1}{2} \ \left( \frac{\pi}{4} \right)^2
\end{equation}

\subsubsection{Analysis with a CAS}

Using a Computer Algebra System, we obtain the following values for small even values of the parameter $n$:
\begin{align*}
I^{(1)}_0 &= \frac{\pi ^2}{32}\\
I^{(1)}_2 &= -\frac{\pi ^2}{32}+\frac{\pi}{4}-\frac{\ln 2}{2}\\
I^{(1)}_4 &= \frac{\pi ^2}{32}-\frac{\pi}{6}-\frac 16 +\frac 23 \ln 2\\
I^{(1)}_6 &= -\frac{\pi ^2}{32}+\frac{13\pi}{60}+\frac{13}{60}-\frac{23}{30} \ln 2\\
I^{(1)}_8 &= \frac{\pi ^2}{32}-\frac{19\pi}{105}-\frac{29}{105}+\frac{88}{105} \ln 2\\
I^{(1)}_{10} &=-\frac{\pi ^2}{32} +\frac{263\pi}{1260}+\frac{2333}{7560}-\frac{563}{630} \ln 2\\
I^{(1)}_{12} &= \frac{\pi ^2}{32}-\frac{1289\pi}{6930}-\frac {3578}{10395} +\frac {3254}{3465} \ln 2\\
I^{(1)}_{14} &= -\frac{\pi ^2}{32} +\frac{36979\pi}{180180}+\frac{397753}{1081080}-\frac{88069}{90090} \ln 2
\end{align*}

\subsection{The case of odd $n$}

\subsubsection{Analytic work}
For a given integer $p \geq \ 1$ and odd $n$, we have:
\begin{equation*}
I^{(p)}_n + 2 \left( I^{(p)}_{n-2}+I^{(p)}_{n-4}+ \dots + I^{(p)}_3 \right) +I^{(p)}_1 =\left( \frac{\pi}{4} \right)^p \underset{k=2}{\overset{n-1}{\sum}} \frac {1}{k}
+ p \underset{k=2}{\overset{n-1}{\sum}}  \frac {I_{k}}{k}
\end{equation*}
By addition and substraction of the iterative formula (i.e., Equation (\ref{fundamentalnp})), we obtain:
\begin{equation*}
-I^{(p)}_n + I^{(p)}_1 = \left( \frac{\pi}{4} \right)^p \underset{k=2}{\overset{n-1}{\sum}} \frac {(-1)^k}{k}
+ p \underset{k=2}{\overset{n-1}{\sum}} \frac {(-1)^k I^{(p-1)}_{k}}{k}
\end{equation*}

For $p=0$, we find the previous expression.
For even $n$ and  $p=1$, we obtain:
\begin{equation*}
I^{(1)}_n + 2 (I^{(1)}_{n-2}+I^{(1)}_{n-4}+ \dots + I^{(1)}_3 ) +I^{(1)}_1 =\frac{\pi}{4}\underset{k=2}{\overset{n-1}{\sum}} \frac {1}{k}
+ \underset{k=2}{\overset{n-1}{\sum}}  \frac {I_k}{k}
\end{equation*}
Then we derive a general expression for $I^{(1)}_n$ and $n=2l+1$:
\begin{proposition}
\label{gen expression odd index}
\begin{equation*}
I^{(1)}_{2l+1} = \frac{\pi}{4} \underset{k=2}{\overset{2l}{\sum}} \frac {(-1)^k}{k}
+  \underset{k=2}{\overset{2l}{\sum}} \frac {(-1)^k I_k}{k}- I_1
\end{equation*}
\end{proposition}

\subsubsection{Analysis with a CAS}
Using a Computer Algebra System, we obtain the following values for small odd values of the parameter $n$:
\begin{align*}
I^{(1)}_1 &= \frac{\pi}{8} \ln 2 + \int_0^{{\pi}{4}} \ln (\cos x)) \; dx\\
I^{(1)}_3 &= -\frac{\pi}{8} \ln 2 - \int_0^{{\pi}{4}} \ln (\cos x)) \; dx + \frac{\pi}{4} -\frac 12\\
I^{(1)}_5 &= \frac{\pi}{8} \ln 2 + \int_0^{{\pi}{4}} \ln (\cos x)) \; dx - \frac{\pi}{4} +\frac 23\\
I^{(1)}_7 &= -\frac{\pi}{8} \ln 2 - \int_0^{{\pi}{4}} \ln (\cos x)) \; dx + \frac{\pi}{3} -\frac {73}{90}\\
I^{(1)}_9 &= \frac{\pi}{8} \ln 2 + \int_0^{{\pi}{4}} \ln (\cos x)) \; dx - \frac{\pi}{3} + \frac{284}{315}\\
I^{(1)}_{11} &= -\frac{\pi}{8} \ln 2 - \int_0^{{\pi}{4}} \ln (\cos x)) \; dx + \frac{23 \pi}{60} -\frac{3103}{3150}\\
I^{(1)}_{13} &= \frac{\pi}{8} \ln 2 + \int_0^{{\pi}{4}} \ln (\cos x)) \; dx - \frac{23 \pi}{60} +\frac{54422}{51975}\\
I^{(1)}_{15} &= -\frac{\pi}{8} \ln 2 - \int_0^{{\pi}{4}} \ln (\cos x)) \; dx + \frac{44 \pi}{105} - \frac{10459489}{9459450}\\
I^{(1)}_{16} &=\frac{\pi}{8} \ln 2 + \int_0^{{\pi}{4}} \ln (\cos x)) \; dx -\frac{44 \pi}{105} + \frac{5452712}{4729725}
\end{align*}
The following remarks can be made:
\begin{enumerate}
\item Two subsequences appear, for even indices and for odd indices.
\item The integral $\frac{\pi}{8} \ln 2 + \int_0^{{\pi}{4}} \ln (\cos x)) \; dx$ appears for an index $n \equiv 1 \mod 4$, and its opposite for
 $n \equiv 3 \mod 4$. The sign change comes from the $(-1)^k$ factor in the two first terms in Proposition (\ref{gen expression odd index}).
\item The integral $\int_0^{{\pi}{4}} \ln (\cos x) \; dx$ is left in closed form and cannot be evaluated analytically.
 \end{enumerate}

 With Maple, the following result is obtained:
 \begin{equation}
 \label{ln cos integral}
 \int_0^{\pi /4} \ln (\cos x)\; dx = \frac{\pi}{4} \ln 2 + \frac 12 G,
 \end{equation}
 where $G$ denotes the Catalan constant\footnote{See the sequence \url{https://oeis.org/A006752}}.

In \cite{adamchik}, the integral presentations of the Catalan constant are given, but without analytic proofs. There, computations have been performed with Mathematica. Note that this is a black-box usage of the CAS, in contradiction to Buchberger's point of view. Nevertheless, we can follow a path in reversed direction and note that
\begin{equation}
G=\underset{r=0}{\overset{\infty}{\sum}} \frac {(-1)^r}{(2r+1)^2}.
 \end{equation}

\section{Study of the influence of  the parameter $p$ }

In this section we investigate the dependence of the parameter $p$ on  the 2-parameter  integral $I_{n}^{(p)}$. The fundamental iterative formula (\ref{fundamentalnp}) connects $I_{n}^{(p)}$ with
the terms $I_{n-2}^{(p)}$ and $I_{n-1}^{(p-1)}$ and this for a given $p$.

The influence of $p$ may be checked in two ways.

 First, in the given integration interval, i.e. for  $0 \ \leq x \leq \ \frac{\pi}{4} \ < 1$,  we have $0<x^p <1$. It follows that $I_{n}^{(p)} < I_n$, where
\begin{equation}
I_n=\int_0^{\pi /4}  \tan^n x\; dx.
\end{equation}

\subsection{The tangent power integral}

This last equation has been studied by the authors in reference (\cite{wallis-catalan}).
In this reference , we derived closed combinatorial formulas for the tangent-power integral , according to the congruence class modulo 4 of the parameter $m$. We presented three different forms of the definite integral:
\begin{enumerate}
\item  \underline{A finite series formula:} for even indices, the following holds:
\begin{equation}
\label{eq 3}
I_{2k} = (-1)^k \frac{\pi}{4} + \underset{l=1}{\overset{k}{\sum}}\frac{(-1)^{(l+k)}}{2l-1}.
\end{equation}
For odd $n$, we have:
\begin{equation*}
 I_{2k+1} = (-1)^k \frac{\ln 2}{2} + \underset{l=1}{\overset{k}{\sum}}\frac{(-1)^{(l+k)}}{2l}
 \end{equation*}
\item \underline{Formulas involving the \emph{Digamma function}:} For even $n$, we have:
\begin{equation}
 I_{2k}  = \Psi \left( k+ \frac 12 \right) -\Psi(2k)+  (-1)^k \frac{\pi}{4}.
\end{equation}
and for odd $n$, the following holds:
\begin{equation*}
I_{2k+1}   =1-\Psi \left( \frac n2 + \frac 12 \right)+\Psi(n)+(-1)^k \log \left( \frac{1}{\sqrt{2}} \right).
\end{equation*}
\item An equivalent formula using factorials:
\begin{align*}
 I_{4k-2} & = - \frac{\pi}{4} +\frac{k}{(2k-1)!!} \\
I_{4k-1} & =   \frac{\ln 2}{2} + \frac{1}{2} +\frac{k}{(2k-1)!!} \\
 I_{4k} & =  \frac{\pi}{4} +\frac{k}{(2k-1)!!} \\
I_{4k+1} & = \frac{\ln 2}{2} + \frac{1}{2} +\frac{k}{(2k)!!}
\end{align*}
\end{enumerate}
These different expressions have different physical meanings leading to different types of understandings and teaching:

\subsection{An iterative algorithm for $I_{n}^{(p)}$ for $n$ and $p=1,...,k$}
In the previous section we derived a general iterative equation for fixed $n$ between $I_{n}^{(p)}$ and the terms
$I_{k}^{(p-1)}, \; k < n$:
\begin{itemize}
\item If $n$ is even:
\begin{equation*}
I^{(p)}_n = \left( \frac{\pi}{4} \right)^p \underset{k=1}{\overset{n-1}{\sum}} \frac {(-1)^k}{k}
+ p \underset{k=1}{\overset{n-1}{\sum}} \frac {(-1)^k I^{(p-1)}_k}{k}- \frac{1}{p+1} \ \left( \frac{\pi}{4} \right)^{(p+1)}
\end{equation*}
\item If $n$ is odd:
\begin{equation*}
-I^{(p)}_n + I^{(p)}_1 = \left( \frac{\pi}{4} \right)^p \underset{k=2}{\overset{n-1}{\sum}} \frac {(-1)^k}{k}
+ p \underset{k=2}{\overset{n-1}{\sum}} \frac {(-1)^k I^{(p-1)}_{k}}{k}
\end{equation*}
\end{itemize}
As a consequence, we may compute iteratively any integral $I_{n}^{(p)}$.
The method consists in:
\begin{enumerate}[label=(\roman*)]
\item Check if $n$ is even or odd
\item If $n$ is even, compute $I_{k}^{(0)}=\int_0^{\pi /4}  \tan^n x\; dx $ for $k < n$.
\item Then use the above formula to compute $I_{n}^{(1)}$
\item Then compute $I_{k}^{(1)}$  for $k < n$.
\item Compute $I_{n}^{(2)}$  and so on until $I_{n}^{(p)}$
\item If $n$ is odd, compute $I_{k}^{(0)}=\int_0^{\pi /4}  \tan^n x\; dx $ for $1< k < n$.
\item Then use the above formula to compute $I_{n}^{(1)}$
\item Then compute $I_{k}^{(1)}$  for $k < n$.
\item Compute $I_{n}^{(2)}$  and so on ....until $I_{n}^{(p)}$
\end{enumerate}

\subsection{Analysis of an integral using a power series}

We consider the integral
\begin{equation}
J_{n}= \int_0^{\pi /4} \frac{ \tan^n x}{1-x} dx.
\end{equation}
For any integer $n>0$, we were unable to obtain a direct analytic computation of this integral using a CAS. Only numerical values have been obtained. Nevertheless, the integral $J_n$ may be computed using the integral $ I_{n}^{(p)}$ studied earlier. The domain of integration is the interval $\left[ 0, \frac{\pi}{4}\right]$, therefore $|x| < 1$ and we have
\begin{equation*}
\frac{1}{1-x} = \sum_{i=0}^{\infty} x^i,
\end{equation*}
and it follows that
\begin{equation*}
J_{n}= \int_0^{\pi /4} \frac{\tan^n x}{1-x}; dx = \sum{i=0}^{\infty} \int_0^{\pi /4} x^i \; \tan^n x \; dx=\sum_{i=0}^{\infty} I_{n}^{i}.
\end{equation*}

The computation of $J_{n}$ requires the computation of the integrals $I_{n}^{i}$ for all $i$. In the previous section, we presented an iterative method to compute $I_{n}^{i}$ from the different $I_{n}^{i-1}; i \geq 1$. Now we present an approximate method to compute $J_{n}$.

We have:
\begin{equation*}
 \forall n \geq n_0, \;  I_{n}^{(i)} \approx (\frac{\pi}{4})^i) \int_0^{\frac{\pi}{4}}  \tan^n( x)\; dx.
\end{equation*}
It follows that:
\begin{equation*}
J_{n}=  \sum_{i=0}^{\infty} I_{n}^{(i)} = \sum_{i=0}^{n_0} I_{n}^{i}+ \sum_{i=n_0}^{\infty} I_{n}^{(i)}.
\end{equation*}
Using the above approximation, we obtain:
\begin{equation*}
J_{n}=   \sum_{i=0}^{n_0} I_{n}^{(i)}+ \left[\sum_{i=n_0}^{\infty} (\frac{\pi}{4})^i)\right] \left[ \int_0^{\frac{\pi}{4}}  \tan^n( x)\; dx \right].
\end{equation*}
Therefore, we obtain:
\begin{equation*}
J_{n}=   \sum_{i=0}^{n_0} I_{n}^{i}+ \left[ (\frac{\pi}{4})^{n_0} \frac{1}{1- \frac{\pi}{4}}   \right] \; \cdot \; \left[ \int_0^{}  \tan^n( x)\; dx \right].
\end{equation*}
The index $n_0$ may be chosen by using an error level $\epsilon$ and require that $x^{n_0} < \epsilon$.

\subsection{Analysis of the integral $L_{n}= \int_0^{1}  \arctan^n( x)\; dx$ }
Note first that this integral is a particular case of the integral $ I_{n}^{(p)}$ studied above. When asked to compute this integral for small values of the parameter, a CAS returns only te closd from of the integral. Only in numerical form, some other output can be obtained, but this is irrelevant to our purpose here..

Consider the change of variable: $u = \arctan(x)$; then $x = \tan(u)$ for $ 0 \le u \le \frac{\pi}{4}$.
Note that $dx = (1+ \tan^2 u) \; du$. With this change of variable $L_{n} = \int_0^{\pi /4} \ u^n \ (1+ \tan^2 u) \; du$, whence:
\begin{equation}
\label{eq Ln}
L_{n} = \int_0^{\pi /4}  u^n \; du + I_{2}^{(p)}
\end{equation}

Finally, we have the following proposition:

\begin{proposition}
For every non negative integers $n$ and $p$, the following holds:
\begin{equation*}
L_{n} = \frac{1}{n+1} \left( \frac{\pi}{4} \right)^{n+1} + I_{2}^{(p)}.
\end{equation*}
\end{proposition}

\section{Some thoughts in the aftermath}

Previous explorations of parametric definite integrals, no matter how many parameters were involved, opened for us various ways to build bridges between mathematical domains. In many cases, integral presentations of combinatorial objects have been derived. In other cases, they provided an opportunity to prove combinatorial identities \cite{stirling,wallis-catalan,Qi et al 2017,Qi Elemente 2017}. A CAS is an important tool for such an exploration and proofs. Nevertheless, it is sometimes hard to conjecture such an identity and another technology can be useful. This was the case for the mentioned works, and using the Online Encyclopedia of Integer Sequences revealed efficient. Of course, it helped to have a conjecture, but this conjecture had to be proven by more theoretical means. 

For the examples of the present work, the situation is somehow different. We could try various inputs for searching the database, but none provided a ready-to-use answer. For all our trials, the database proposed a modified version of the input, and understanding the proposed sequences was non-trivial.  This is a good illustration that exploring mathematical objects in a technology-rich environment does not follow a predefined path.

The OECD defined the so-called 4 C's of 21st century education: Communication, Collaboration, Critical Thinking and Creativity. These are crucial in various settings (see for example \cite{saimon et al}. Here, the collaboration can be between humans, but this is not central. Collaboration with technology, i.e. human and machine together, requires strong Critical Thinking in order to analyse the outputs of the CAS and of the interactive database.  Critical Thinking and Creativity are central skills, opposite to Black-Box usage of technology. Further, Creativity leads to the possible conjectures. This is a modern way to deal with mathematics, with a scheme exploration-conjecture-proof. It request also the understanding that the new technological skills are an integral part of the new mathematical knowledge, and that a new technological discourse has to be developed \cite{artigue,MDZ}. Such exploration can be initiated by a 5th C, namely Curiosity. This has been the basis of \cite{soil}, where mathematical curiosity and communication and collaboration between researchers yielded an interesting application.
   
Nowadays, curiosity can push in new  directions. We asked an AI for references about parametric definite integrals. Most answers gave general chapters in Calculus books, the only relevant to \emph{parametric} definite integrals was to present authors' presentation in a conference \cite{DPZ-ACA 2019}.

Finally, we wish to recall what we wrote in Section \ref{intro}. Parametric definite integrals enable students to step forwards towards more abstract understanding than what is offered in a first Calculus course. But they are not always a pure mathematical exercise.  They have numerous applications in science, in particular in Physics; \cite{Hu and Rebello - Physics} provides an insight into students' difficulties to translate the knowledge acquired in definite integrals in a Calculus course into the needed skills for Physics. Nevertheless, the topics is quite rare in Calculus textbooks; we found a subchapter in \cite{Lax and Terell} (subchapter 7.4). Regarding the mental schemes that they require, these are parallel to what is required in Linear Algebra for solving parametric systems of linear equations; see \cite{matricial computations}, where some problems are analyzed, such as the lack of investigation of special values of the parameters\footnote{Such problems, when working with AI for mathematics education, are analyzed in \cite{bagno et al}.}.
\vskip 0.2cm

{\bf Acknowledgements:}
The 2nd author has been partially supported by the CEMJ Chair at JCT.

{\bf Declaration:}
The authors declare no conflict of interest.

\bibliographystyle{amsalpha}

\normalsize
\end{document}